\documentclass[12pt]{article}

\usepackage{amssymb}
\usepackage{amsmath}

\newcommand{\lap}{\mbox{$\bigtriangleup$}}

\newcommand{\ra}{{\mbox{$\rightarrow$}}}
\newcommand{\be}{\begin{equation}}
\newcommand{\ee}{\end{equation}}

\newcommand{\F}{{\cal F}}
\newcommand{\sch}{{\cal S}}

\newcommand{\mint}{-\!\!\!\!\!\!\!\int}

\newtheorem{theorem}{Theorem}[section]

\newtheorem{mrem}[theorem]{Remark}

\newtheorem{mthm}[theorem]{Theorem}
\newtheorem{mpro}[theorem]{Proposition}

\newtheorem{thm}[theorem]{Theorem}
\newtheorem{pro}[theorem]{Proposition}
\newtheorem{lem}[theorem]{Lemma}
\newtheorem{cor}[theorem]{Corollary}

\begin{document}

\title{Some Liouville theorems for \\ the fractional Laplacian}
\medskip

\author{
Wenxiong Chen\thanks{Partially supported by the Simons Foundation Collaboration Grant for Mathematicians 245486.}, \hskip .1in Lorenzo D'Ambrosio\thanks{Partially supported by the Italian MIUR National Research Project: Quasilinear Elliptic
Problems and Related Questions} \hskip .1in  and  \hskip .1in  Yan Li}

\date{\today}
\maketitle
\rightline{{\em Dedicated to Professor Enzo Mitidieri}}
\rightline{{\em on the occasion of his 60th birthday }}

\begin{abstract}
In this paper, we prove the following result.
Let $\alpha$ be any real number between $0$ and $2$.  Assume that $u$ is a solution of
$$ \left\{\begin{array}{ll}
(-\lap)^{\alpha/2} u(x) = 0 , \;\; x \in \mathbb{R}^n ,\\
\displaystyle\underset{|x| \ra \infty}{\underline{\lim}} \frac{u(x)}{|x|^{\gamma}} \geq 0 ,
\end{array}
\right.
$$
for some $0 \leq \gamma \leq 1$ and $\gamma < \alpha$.  Then $u$ must be constant throughout $\mathbb{R}^n$.

This is a Liouville Theorem for $\alpha$-harmonic functions under a much weaker condition.

For this theorem we have two different proofs by using two different methods:
One is a direct approach using potential theory.
The other is by Fourier analysis as a corollary of the fact that
the only $\alpha$-harmonic  functions  are affine.
\end{abstract}

\bigskip


\section{Introduction}

The well-known classical Liouville's Theorem states that

{\em Any harmonic function bounded below in all of $R^n$ is constant.}

One of its important applications is the proof of the Fundamental Theorem of Algebra. It is
also a key ingredient in deriving a priori estimates for solutions in PDE analysis.

This Liouville Theorem has been generalized to the fractional Laplacian by Bogdan, Kulczycki, and Nowak \cite{BKN}:
\begin{mpro} Let $0<\alpha<2$ and $n \geq 2$. Assume that $u$ is a solution of
$$\left\{\begin{array}{ll}
(-\lap)^{\alpha/2} u(x) = 0 , & x \in \mathbb{R}^n , \\
u(x) \geq 0 , & x \in \mathbb{R}^n .
\end{array}
\right.
$$
Then $u$ must be constant.
\end{mpro}

The same result has been proved by Zhuo, Chen, Cui, and Yuan \cite{ZCCY} using a completely different method;
and then interesting applications of this Liouville theorem to integral representations of solutions for nonlinear equations and systems involving the fractional Laplacian were investigated in the same article.

The above proposition has also been proved but not explicited stated in \cite{L} for $\alpha$-harmonic functions in the average sense (see also Section \ref{appendix} for the definition).
Indeed it can be deduced from \cite[Theorem 1.30]{L}.

\medskip

In \cite{ABR}, Axler,
 Bourdon, and Ramey replaced the condition ``bounded below'' in the classical Liouville Theorem by a much weaker one:
$$\underset{|x| \ra \infty}{\underline{\lim}} \frac{u(x)}{|x|} \geq 0 .$$

Then Enzo Mitidieri \cite{Mi1} conjectured that a similar result should hold for the fractional Laplacian.
The main purpose of this paper is to prove this conjecture.

The fractional Laplacian in $R^n$ is a nonlocal pseudo-differential operator, taking the form
\begin{equation}
(-\Delta)^{\alpha/2} u(x) = C_{n,\alpha} \, \lim_{\epsilon \ra 0} \int_{R^n\setminus B_{\epsilon}(x)} \frac{u(x)-u(z)}{|x-z|^{n+\alpha}} dz,
\label{Ad7}
\end{equation}
 where $\alpha$ is any real number between $0$ and $2$.
This operator is well defined in $\sch$, the Schwartz space of rapidly decreasing $C^{\infty}$
functions in $R^n$. In this space, it can also be equivalently defined in terms of the Fourier transform
$$ \F({(-\Delta)^{\alpha/2} u}) (\xi) = |\xi|^{\alpha} \F(u) (\xi), $$
where $\F(u)$ is the Fourier transform of $u$. One can extend this operator to a wider space of distributions.

Let
$$L_{\alpha}=\{u: R^n\rightarrow R \mid \int_{R^n}\frac{|u(x)|}{1+|x|^{n+\alpha}} \, d x <\infty\}.$$

Then in this space, one can defined $(-\Delta)^{\alpha/2} u$ as a distribution by
$$< (-\Delta)^{\alpha/2}u(x), \phi> \, = \, \int_{R^n}  u(x) (-\Delta)^{\alpha/2} \phi (x) dx ,  \;\;\; \forall \, \phi \in C_0^{\infty}(R^n) . $$

The operator can be also defined by considering the following problem in $\mathbb R^n\times[0,+\infty)$:
$$ f(x,0)=u(x)\ on\ \mathbb R^n,\qquad \lap_xf+\frac{1-\alpha}{y}f_y+f_{yy}=0
$$
and defining
$$C{(-\Delta)^{\alpha/2} u} (x):=-\lim_{y\to 0} y^{1-\alpha}f_y(x,y),$$
where $C=C(n,\alpha)$ is a suitable positive constant. See \cite{CS} for more details.

\medskip

Throughout this paper, we will consider the fractional Laplacian defined in (\ref{Ad7}). We say that $u$ is $\alpha$-harmonic if $u \in L_{\alpha}$ and satisfies
$$(- \Delta)^{\alpha/2} u = 0$$
in the sense of distributions:
$$ \int_{R^n}  u(x) (-\Delta)^{\alpha/2} \phi (x) dx = 0 \;\; \forall \, \phi \in C_0^{\infty}(R^n) .$$

Our main objective is to prove the following result.

\begin{mthm} Let $0<\alpha<2$. Assume that  $u \in L_{\alpha}$ is $\alpha$-harmonic and
\begin{equation}
\displaystyle\underset{|x| \ra \infty}{\underline{\lim}}\frac{u(x)}{\mid x\mid^{\gamma}}\geq0,
\label{a2}
\end{equation}
for some $0\leq\gamma\leq1$ and $\gamma<\alpha$.

 Then $$u(x) \equiv C, \;\; x \in \mathbb{R}^{n}.$$
 \label{mthm1}
\end{mthm}

We prove the above result by two different methods.
The first one is given in Section \ref{2} and it is a direct proof based on classical potential theory.

The second proof is a direct consequence of the following
\begin{mthm}\label{th:affine} Let $0<\alpha<2$, and $u\in L_\alpha$ be $\alpha$-harmonic.
Then $u$ is affine. In particular, if $\alpha\le 1$, then $u$ is constant.
\end{mthm}
The proof of Theorem \ref{th:affine} is based on Fourier analysis and Section \ref{3} is devoted to its proof.
Notice that the tool of Fourier transform is not always available in different frameworks. This is the main reason why we present an alternative proof based on
potential theory which, in turn, is well developed even in the metric space framework.

We notice that the first proof is more articulate than the second one.
It allows us to develop some tools which are of interest in themselves. Some key ingredients used in this proof are more or less well-known, however, for readers' convenience, we provide details in the Section \ref{appendix}. There, we will also derive that (see Corollary \ref{corA}) if $u \in L_{\alpha}$ is $\alpha$-harmonic in the sense of distributions, then it is smooth, and it is also $\alpha$-harmonic in the sense of average. The converse can be deduced from the
Pizzetti type formula for the fractional Laplacian (see Lemma \ref{lem3.2}). This formula was proved in \cite{L}, here we provide a different and more direct argument.
We point out  the importance of this kind of results
in the study of polyharmonic operator
(see \cite{P1,P2,Sb} and more recently \cite{Boj}),
 sublaplacian on Carnot groups (see \cite{Bon}) and
the heat operator (see \cite{DR}).

For more related articles on the uniqueness or non-existence of solutions for nonlinear elliptic equations-- Liouville type theorems--and their applications, in particular, for those non-local equations, please see \cite{CDM, CLO, CLO1, DS, FC, LiZ, M, MC, Mi, Z} and the references therein.

Finally, we notice that after a partial version of these paper circulated \cite{Fa} have proved
 Theorem \ref{th:affine} with different methods than ours.

\medskip

\noindent{\bf Note.} On May 4th, Professor Enzo Mitidieri gave this problem to us through a private communication. A couple of weeks later, we have proved this conjecture. The results and the ideas of proofs have been presented by the first author at the Special Session on Nonlinear Elliptic Partial Differential Equations and Systems at The 10th AIMS Conference on DSDE on July 9th in Madrid. Now this paper has been
accepted for publication in Nonlinear Analysis Analysis: Theory, Methods \& Applications.
\medskip

The authors are very grateful to Professor Enzo Mitidieri for suggesting this interesting problem to us and for his insightful discussions.

\section{Proof of Theorem \ref{mthm1}}\label{2}


{\em Step 1.}

We first show that, for $|x|<r$,
 \begin{equation}
u(x)= \underset{|y|>r}{\int}P_{r}(y,x)u(y)dy,
\label{P}
\end{equation}
where $P_r(y,x)$ is some kind of Poisson kernel for $|x|<r$:
\begin{equation}
P_{r}(y,x)=\left\{\begin{array}{ll}
\frac{\Gamma(n/2)}{\pi^{\frac{n}{2}+1}} \sin\frac{\pi\alpha}{2}
\left[\frac{r^{2}-|x|^{2}}{|y|^{2}-r^{2}}\right]^{\frac{\alpha}{2}}\frac{1}{|x-y|^{n}},& |y|>r,\\
0,& |y|<r.
\end{array}
\right.
\end{equation}

Let \begin{equation}
\hat{u}(x)=\left\{\begin{array}{ll}
\int_{|y|>r}P_{r}(y,x)u(y)dy,& |x|<r,\\
u(x),& |x|\geq r.
\end{array}
\right.
\label{z6}
\end{equation}
Then we have

\begin{pro}

$\hat{u}$ is $\alpha$-harmonic in the ball $B_r(0)$.
\label{pro2.1}
\end{pro}

In \cite{L}, Landkof pointed out that, by using a similar method in his book, one can prove the result of this proposition. Since the proof was not given there, and it is quite long and complex, for reader's convenience, we will
present it in the Section \ref{appendix}.

Let $w(x)=\hat{u}(x)-u(x)$, then by Proposition \ref{pro2.1}, we derive
\begin{equation}
\left\{\begin{array}{ll}
(-\lap)^{\frac{\alpha}{2}}w(x)=0, &|x| < r,\\
w(x)=0, &|x| \geq r.
\end{array}
\right.
\end{equation}

To show that $w \equiv 0$, we employ the following {\em maximum principle}.

\begin{lem}
(Silvestre \cite{Si})
 \quad Let $\Omega$ be a bounded domain in $\mathbb{R}^{n}$, and assume that $v$ is a lower semi-continuous function on $\overline{\Omega}$ satisfying
\begin{equation}
\left\{\begin{array}{ll}
(-\lap)^{\frac{\alpha}{2}}v\geq0  &\mbox{in } \Omega,\\
v\geq0&\mbox{on } \mathbb{R}^{n}\backslash\Omega.
\end{array}
\right.
\end{equation}
then $v\geq0$ in $\Omega$.
\label{max}
\end{lem}

Applying this lemma to both $v=w$ and $v=-w$, we conclude that
 $$w(x)\equiv 0.$$
Hence $$\hat{u}(x)\equiv u(x).$$
This verifies (\ref{P}).

{\em Step 2.}

We will show that, for each fixed $x\in\mathbb{R}^{n}$, and for any unit vector ${\bf \nu}$, we have
\begin{equation}
\frac{\partial u}{\partial\nu}(x)\geq0.
\label{b1}
\end{equation}

It follows from the arbitrariness of $\nu$ that
$$\nabla u(x)=0,$$
and therefore
$$u\equiv C \mbox{ in } \mathbb{R}^{n}.$$

Now what left is to prove (\ref{b1}).  Through an elementary calculation, one can derive that
 \begin{equation}
\frac{\partial u}{\partial x_{i}}(x)
=-\int_{|y|>r}P_{r}(y,x)\left[\frac{\alpha x_{i}}{r^{2}-|x|^{2}}+\frac{n(x_{i}-y_{i})}{|y-x|^{2}}\right]u(y)dy,
\end{equation}
and consequently,
 \begin{equation}
\frac{\partial u}{\partial\nu}(x)
=-\int_{|y|>r}P_{r}(y,x)\left[\frac{\alpha x\cdot\nu}{r^{2}-|x|^{2}}+\frac{n(x-y)\cdot\nu}{|y-x|^{2}}\right]u(y)dy.
\end{equation}

By (\ref{a2}), we see that for any fixed $\epsilon>0$, when $|x|$ is sufficiently large, it holds
\begin{equation}
u(x)\geq -\epsilon|x|^{\gamma}.
\label{a12}
\end{equation}
Otherwise, there exists an $\epsilon_{o}$ and a sequence $\{x^{k}\}$ with $|x^{k}|\rightarrow\infty$ as $k\rightarrow\infty$ such that
\begin{equation}
u(x^{k})< -\epsilon_{o}|x^{k}|^{\gamma},
\end{equation}
and hence
\begin{equation}
\lim_{k \ra \infty} \frac{u(x^{k})}{\mid x^{k}\mid^{\gamma}}<-\epsilon_{o}.
\end{equation}
This is a contradiction with our assumption (\ref{a2}).

Now write
\begin{eqnarray}
\frac{\partial u}{\partial\nu}(x)
&=&-\int_{|y|>r}P_{r}(y,x)\left[\frac{\alpha x\cdot\nu}{r^{2}-|x|^{2}}+\frac{n(x-y)\cdot\nu}{|y-x|^{2}}\right][u(y)+\epsilon|y|^{\gamma}]dy\nonumber\\
&+&\int_{|y|>r}P_{r}(y,x)\left[\frac{\alpha x\cdot\nu}{r^{2}-|x|^{2}}+\frac{n(x-y)\cdot\nu}{|y-x|^{2}}\right]\epsilon|y|^{\gamma}dy\nonumber\\
&=&H_{1}+H_{2}.\label{b6}
\end{eqnarray}

For each fixed $x$, let $r$ be sufficiently large to ensure both (\ref{a12}) and (\ref{P}) and the following
$$\left|\frac{\alpha x\cdot\nu}{r^{2}-|x|^{2}}\right|\leq\frac{1}{r}; \;\; \mbox{and for } |y|>r, \;\; \left| \frac{n(x-y)\cdot\nu}{|y-x|^{2}}\right|\leq\frac{n}{|y-x|}\leq\frac{2n}{r} .$$

Then we have
\begin{eqnarray}
H_{1}&\geq&-\frac{2n+1}{r}\int_{|y|>r}P_{r}(y,x)[u(y)+\epsilon|y|^{\gamma}]dy\nonumber\\
&=&-\frac{2n+1}{r}u(x)-\frac{2n+1}{r}\int_{|y|>r}P_{r}(y,x)\epsilon|y|^{\gamma}dy\nonumber\\
&=&-\frac{2n+1}{r}u(x)-H_{11}.\label{b8}
\end{eqnarray}

Obviously, the first term on the right hand side of the above inequality approaches $0$ as $r \rightarrow \infty$. Next we show that one also has
\begin{equation}
H_{11}<C\epsilon, \;\; \mbox{ for all sufficiently large $r$. }\label{b9}
\end{equation}
Here and below, the letter $C$ stands for various constants.
In fact,
\begin{eqnarray}
H_{11} &=&\frac{2n+1}{r}\int_{|y|>r}P_{r}(y,x)\epsilon|y|^{\gamma}dy \nonumber\\
&=&\frac{2n+1}{r} \, \frac{\Gamma(n/2)}{\pi^{\frac{n}{2}+1}} \sin\frac{\pi\alpha}{2}(r^{2}-|x|^{2})^{\frac{\alpha}{2}}\epsilon
\int_{|y|>r}\frac{|y|^{\gamma}}{(|y|^{2}-r^{2})^{\frac{\alpha}{2}}|x-y|^{n}}dy \nonumber\\
&=&C\epsilon\frac{(r^{2}-|x|^{2})^{\frac{\alpha}{2}}}{r}
\int_{|y|>r}\frac{|y|^{\gamma}}{(|y|^{2}-r^{2})^{\frac{\alpha}{2}}|x-y|^{n}}dy \nonumber\\
&\leq&C\epsilon\frac{(r^{2}-|x|^{2})^{\frac{\alpha}{2}}}{r}
\int_{|y|>r}\frac{|y|^{\gamma}}{(|y|^{2}-r^{2})^{\frac{\alpha}{2}}(|y|-|x|)^{n}}dy \nonumber\\
&=&C\epsilon\frac{(r^{2}-|x|^{2})^{\frac{\alpha}{2}}}{r}
\int_{r}^{\infty}{\frac{\tau^{\gamma+n-1}}{(\tau^{2}-r^{2})^{\frac{\alpha}{2}}(\tau-|x|)^{n}}}d\tau\label{b3}\\
&=&C\epsilon\frac{(r^{2}-|x|^{2})^{\frac{\alpha}{2}}}{r^{\alpha-\gamma+1}}
\int_{1}^{\infty}{\frac{s^{\gamma+n-1}}{(s^{2}-1)^{\frac{\alpha}{2}}(s-\frac{|x|}{r})^{n}}}ds\label{b4}\\
&\leq&  \frac{C\epsilon}{r^{1-\gamma}} \int_1^\infty \frac{s^{\gamma +n-1}}{(s^2-1)^{\frac{\alpha}{2}}
 (s -\frac{|x|}{r} )^n} d s  \label{b5}
\end{eqnarray}
  We derive (\ref{b3}) and (\ref{b4}) by letting $|y|=\tau$ and $\tau=rs$ respectively. From the assumption that $\gamma<\alpha$ and $0<\alpha<2$,  it is not difficult to see that
 \begin{equation}
 \int_1^\infty \frac{s^{\gamma +n-1}}{(s^2-1)^{\frac{\alpha}{2}}
 (s -\frac{|x|}{r} )^n} d s <\infty.
 \label{1.1}
 \end{equation}
 Noticing that $\gamma \leq 1$, and by (\ref{b5}), we arrive at (\ref{b9}).

Through an identical argument, we can deduce
$$ |H_2| \leq  \frac{2n+1}{r}\int_{|y|>r}P_{r}(y,x)\epsilon|y|^{\gamma}dy ,$$
and
\begin{equation}
|H_2| < C \epsilon , \;\; \mbox{for sufficiently large r} .
\label{1.2}
\end{equation}

From (\ref{b6}), (\ref{b8}),  (\ref{b9}), and (\ref{1.2}), letting $r \ra \infty$ and we conclude that
$$\frac{\partial u}{\partial\nu}\geq C\epsilon.$$
The fact that $\epsilon$ is arbitrary establishes (\ref{b1}), hence proves the theorem.

\section{Proof of Theorem \ref{th:affine} and remarks.}\label{3}







In this section we shall prove Theorem \ref{th:affine}.

{\noindent \sl Strategy.} A function $u\in L_\alpha$ is a tempered distribution, and hence $u$ admits a Fourier transform
$\F u$.
If we prove that $\F u$ has support at one point (the origin) then $\F u$ is a finite combination of the Dirac's delta measure and its derivatives. Therefore $u$ is a polynomial.
The only polynomial belongs to $L_\alpha$ are the constants and eventually the affine functions (depending on $\alpha$).

\bigskip
{\noindent \bf Proof.}
From
$$ \langle (-\Delta)^{\alpha/2} u, \psi\rangle = \int_{R^n} u(x)  (-\Delta)^{\alpha/2}\psi(x) dx
  \qquad \forall \psi\in C^\infty_0(R^n),$$
and the fact that
$$\F( (-\Delta)^{\alpha/2} \psi)(\xi) =  |\xi|^\alpha \F(\psi)(\xi)\qquad\quad for\ \ \psi\in \sch,$$
we  observe that $(-\Delta)^{\alpha/2} u =0$ means
 that for any $\psi\in \sch$
\begin{equation} 0 = \langle (-\Delta)^{\alpha/2} u, \psi\rangle = \int_{R^n} u(x)  (-\Delta)^{\alpha/2}\psi(x) dx
=\int_{R^n} u(x)  \F^{-1}(|\xi|^\alpha\F\psi)(x) dx.  \label{eq:1} \end{equation}

We claim that
\begin{equation}\label{eq:3}
   \langle\F u,\phi\rangle =0\qquad for\ \ any\ \  \phi\in \ C_0^\infty(R^n\setminus\{0\}).
\end{equation}
Indeed let $\phi\in \ C_0^\infty(R^n\setminus\{0\})$.
The function ${\phi(\xi)}/{|\xi|^\alpha}$ belongs to $C_0^\infty(R^n\setminus\{0\})\subset \sch$.
Therefore there exists
$\psi\in \sch$ such that $\F(\psi)(\xi) = {\phi(\xi)}/{|\xi|^\alpha}$.

Now, since $u$ is a tempered distribution and from (\ref{eq:1}), we have
$$ \langle\F u,\phi\rangle = \langle\F u,|\xi|^\alpha \F\psi \rangle= $$
$$  = \langle u, \F(|\xi|^\alpha \F\psi )\rangle= \int_{R^n} u(x) \F(|\xi|^\alpha \F\psi )(x)dx=0.$$
That is the claim.

From (\ref{eq:3}) we infer that $\F u$ has support in $\{0\}$.
Then $\F u$ is a finite combination of the Dirac's delta measure at the origin and its derivatives. Therefore, we conclude  that $u$ is a polynomial (see for instance \cite{G}).
Since  $u\in L_\alpha$, we obtain that $u$ has at most a linear growth.
This concludes the proof.

\bigskip
 The same idea can be applied to prove that the polyharmonic tempered distributions are polynomials.

 \begin{mthm} Let $u\in \sch'(R^n)$, $m$ an integer and $(-\lap)^m u =0$.
  Then $u$ is a polynomial.
 \end{mthm}

The proof is similar. Actually the proof is simpler since the symbol of $(-\lap)^m$,
 $|\xi|^{2m}$, being a polynomial, is a good multiplier for tempered distribution. We omit the details.

\begin{mrem} With the same technique one can prove that if $u\in L_\alpha\cap L_\beta$,
with $0<\alpha,\beta<2$ and $u$ solves
$$(-\lap)^{\alpha/2} u(x) +  (-\lap)^{\beta/2} u(x) = 0,\qquad$$
then $u$ is affine (and constant if $\min\{\alpha,\beta\}\le 1$).

A more general result can be obtained with some hypothesis on the symbol of the operator involved. We leave the details to the interested readers.
\end{mrem}

\section{Poisson Representations}\label{appendix}

Let
\begin{equation}
\hat{u}(x)=\left\{\begin{array}{ll}
\int_{|y|>r}P_{r}(y,x)u(y)dy,& |x|<r,\\
u(x),& |x|\geq r,
\end{array}
\right.
\label{hat}
\end{equation}
where
\begin{equation}
P_{r}(y,x)=\left\{\begin{array}{ll}
\frac{\Gamma(n/2)}{\pi^{\frac{n}{2}+1}} \sin\frac{\pi\alpha}{2}
\left[\frac{r^{2}-|x|^{2}}{|y|^{2}-r^{2}}\right]^{\frac{\alpha}{2}}\frac{1}{|x-y|^{n}}, &|y|>r,\\
0, &|y|<r.
\end{array}
\right.
\end{equation}
In this section, we prove
 \begin{thm} Let $u \in L_{\alpha}$ be $\alpha$-harmonic and let
   $\hat{u}(x)$ be defined by (\ref{hat}).
   Then $\hat u$ is $\alpha$-harmonic in $B_r(0)$.
   \label{thm3.1}
 \end{thm}

 The proof consists of two parts. First we show that $\hat{u}$ is harmonic in the average sense (Lemma \ref{lem3.1}), then we show that it is $\alpha$-harmonic by the Pizzetti type formula (Lemma \ref{lem3.2}).

Let
\begin{equation}
\varepsilon^{(r)}_{\alpha}(x)=\left\{\begin{array}{ll}
0, &|x|<r.\\
\frac{\Gamma(n/2)}{\pi^{\frac{n}{2}+1}} \sin\frac{\pi\alpha}{2}\frac{r^\alpha}{(|x|^2-r^2)^{\frac{\alpha}{2}}|x|^n}, &|x|>r.
\end{array}
\right.
\end{equation}

We say that $u$ is $\alpha$-harmonic in the average sense (see \cite{L}) if for small $r$,
$$\varepsilon^{(r)}_{\alpha}\ast u(x)=u(x),$$
where $\ast$ is the convolution.

\begin{lem}
 Assume that
 \begin{equation}
 u \in L_{\alpha}.
 \label{d1}
\end{equation}
 Then $\hat{u}(x)$ is $\alpha$-harmonic in the average sense in $B_{r}(0)$.
\label{lem3.1}
\end{lem}

The following result is the analog of Pizzetti's formula for the classical Laplace operator, namely
 \begin{equation}
\lim_{r\rightarrow0}\frac{1}{r^2}\left[u(x)- \mint_{B_r(x)} u\right]
=c\ (-\lap)u(x).
\label{piz}
\end{equation}
\begin{lem}[Pizzetti's formula] Let $u\in L_\alpha$ and $C^{1,1}$ in a neighborhood of $x$.
Then we have
\begin{equation}
\lim_{r\rightarrow0}\frac{1}{r^\alpha}\left[u(x)-\varepsilon^{(r)}_{\alpha}\ast u(x)\right]
=c\ (-\lap)^\frac{\alpha}{2}u(x).
\label{c0}
\end{equation}
where $c=\frac{\Gamma(n/2)}{\pi^{\frac{n}{2}+1}} \sin\frac{\pi\alpha}{2}$.
\label{lem3.2}
\end{lem}

This formula was obtained in \cite{L}. Here we present a different and more direct proof.
\medskip

Combining the above two Lemmas and the smoothness of $\hat{u}$, we derive
\begin{cor} If $u \in L_{\alpha}$ is $\alpha$-harmonic in the sense of distributions, then it is smooth, and it is also $\alpha$-harmonic in the sense of average. The converse is also true.
\label{corA}
\end{cor}

\noindent\textbf{Proof of Lemma \ref{lem3.1}.}
The outline is as follows.

\textit{i)} \quad Approximate $u$ by a sequence of smooth, compactly supported functions $\{u_k\}$, such that $u_k(x) \ra u(x)$ and
\begin{equation}
\int_{|z|>r}\frac{|u_k(z)-u(z)|}{|z|^n(|z|^2-r^2)^\frac{\alpha}{2}}dz \ra 0.
\label{d2}
\end{equation}
This is possible under our assumption (\ref{d1}).

\textit{ii)} \quad For each $u_k$, find  a signed measure $\nu_k$ such that $supp\,\nu_k\subset B^c_r$ and $$u_k(x)=U_{\alpha}^{\nu_k}(x), \;\; |x|>r.$$
Then $$\hat{u_k}(x)=U_{\alpha}^{\nu_k}(x), \quad |x|<r.$$

\textit{iii)} \quad It is easy to see that $\hat{u_k}(x)$ is $\alpha$-harmonic in the average sense for $|x|<r$.
That is, for each fixed small $\delta>0$,
\begin{equation}
(\varepsilon_\alpha^{(\delta)}\ast\hat{u}_k)(x)=\hat{u}_k(x).
\end{equation}
By showing that as $k \ra \infty$
$$\varepsilon_\alpha^{(\delta)}\ast\hat{u}_k \ra \varepsilon_\alpha^{(\delta)}\ast\hat{u},$$
and
$$\hat{u}_k \ra \hat{u},$$
we arrive at
$$(\varepsilon_\alpha^{(\delta)}\ast\hat{u})(x)=\hat{u}(x),\quad  |x|<r.$$
\bigskip

Now we carry out the details.

\textit{i)} \quad There are several ways to construct such a sequence $\{u_k\}$. One is to use the mollifier.
Let
\begin{equation}
u|_{B_k}(x)=\left\{\begin{array}{ll}
u(x),&|x|<k,\\
0,&|x|\geq k,
\end{array}
\right.
\end{equation}
and
\begin{equation}
J_\epsilon(u|_{B_k})(x)=\int_{\mathbb{R}^n}j_\epsilon(x-y)u|_{B_k}(y)dy.
\end{equation}

For any $\delta>0$, let $k$ be sufficiently large (larger than $r$) such that
\begin{equation}
\int_{|z|\geq k}
\frac{|u(z)|}{|z|^n(|z|^2-r^2)^{\frac{\alpha}{2}}}dz<\frac{\delta}{2}.
\end{equation}
For each such $k$, choose $\epsilon_k$ such that
\begin{equation}
\int_{B_{k+1}\backslash B_r}
\frac{|u_k(z)-u|_{B_k}(z)|}{|z|^n(|z|^2-r^2)^{\frac{\alpha}{2}}}dz
<\frac{\delta}{2},
\end{equation}
where $u_k=J_{\epsilon_k}(u|_{B_k})$. It then follows that
\begin{eqnarray*}
&\quad&\int_{|z|>r}\frac{|u_k(z)-u(z)|}{|z|^n(|z|^2-r^2)^{\frac{\alpha}{2}}}dz\\
&\leq&\int_{B_{k+1}\backslash B_r}\frac{|u_k(z)-u|_{B_k}(z)|+
|u|_{B_k}(z)-u(z)|}{|z|^n(|z|^2-r^2)^{\frac{\alpha}{2}}}dz\\
&&+\int_{|z|>k+1}\frac{|u(z)|}{|z|^n(|z|^2-r^2)^{\frac{\alpha}{2}}}dz\\
&=&\int_{B_{k+1}\backslash B_r}\frac{|u_k(z)-u|_{B_k}(z)|}{|z|^n(|z|^2-r^2)^{\frac{\alpha}{2}}}dz+
\int_{|z|\geq k}\frac{|u(z)|}{|z|^n(|z|^2-r^2)^{\frac{\alpha}{2}}}dz\\
&<&\frac{\delta}{2}+\frac{\delta}{2}=\delta.
\end{eqnarray*}

Therefore, as $k \ra \infty$,
\begin{equation}
\int_{|z|>r}\frac{|u_k(z)-u(z)|}{|z|^n(|z|^2-r^2)^\frac{\alpha}{2}}dz \ra 0.
\label{w1}
\end{equation}

\textit{ii)} \quad For each $u_k$, there exists a signed measure $\psi_k$ such that
\begin{equation}
u_k(x)=U_{\alpha}^{\psi_k}(x).
\end{equation}
 Indeed, let $\psi_k(x)= C (-\lap)^{\frac{\alpha}{2}}u_k(x)$, then
\begin{eqnarray}
U_{\alpha}^{\psi_k}(x)&=&\int_{\mathbb{R}^n}\frac{C}{|x-y|^{n-
 \alpha}}(-\lap)^{\frac{\alpha}{2}}u_k(y)dy\\
 &=&\int_{\mathbb{R}^n}(-\lap)^{\frac{\alpha}{2}}\left[\frac{C}{|x-y|^{n-
 \alpha}}\right]u_k(y)dy\label{d3}\\
 &=&\int_{\mathbb{R}^n}\delta(x-y)u_k(y)dy=u_k(x).
 \end{eqnarray}
Here we have used the fact that $\frac{C}{|x-y|^{n-\alpha}}$ is the fundamental solution of $(-\lap)^{\alpha/2}$.

Let $\psi_k|_{B_r}$ be the restriction of $\psi_k$ on $B_r$ and
\begin{equation}
\tilde{\psi_k}(y)=\int_{|x|<r}P_r(y, x)\psi_k|_{B_r}(x)dx,
\end{equation}
we have
 $$U_{\alpha}^{\tilde{\psi_k}}(x)=U_{\alpha}^{\psi_k|_{B_r}}(x), \;\; |x|> r, $$ and
$supp\,\tilde{\psi_k}\subset B^c_r.$ Here we use the fact (see (1.6.12$'$) \cite{L}) that
\begin{equation} \frac{1}{|z-x|^{n-\alpha}}=\int_{|y|>r}\frac{P_r(y,x)}{|z-y|^{n-\alpha}}dy,
 \:|x|<r,\: |z|>r.
 \label{Po}
 \end{equation}

Let $\nu_k=\psi_k-\psi_k|_{B_r}+\tilde{\psi_k}$, then $supp\,\nu_k\subset B^c_r$, and
$$U_{\alpha}^{\nu_k}(x)=U_{\alpha}^{\psi_k}(x)+U_{\alpha}^{\tilde{\psi_k}}(x)
-U_{\alpha}^{\psi_k|_{B_r}}(x)=U_{\alpha}^{\psi_k}(x),\quad |x|>r.$$
That is
$$ u_k(x) = U_{\alpha}^{\nu_k}(x) ,\quad |x|>r.$$

Again by using (\ref{Po}), we deduce
$$ \hat{u}_k(x) = U_{\alpha}^{\nu_k}(x) , \quad |x|<r.$$
Now it is easy to verify that (see \cite{L}) $\hat{u}_k$ is $\alpha$-harmonic (in the sense of average) in the region $|x|<r$.

\textit{iii)}\quad For each fixed $x$, we first have
$$\hat{u}_k(x) \ra \hat{u}(x).$$
In fact, by (\ref{w1}),
\begin{eqnarray*}
\hat{u}_k(x)-\hat{u}(x)&=&\int_{|y|>r}P_r(y,x)[u_k(y)-u(y)]dy\\
&=&C\int_{|y|>r}\frac{(r^2-|x|^2)^{\frac{\alpha}{2}}[u_k(y)-u(y)]}{(
|y|^2-r^2)^{\frac{\alpha}{2}}|x-y|^n}dy \ra 0.
\end{eqnarray*}
 Next, we show that, for each fixed $\delta>0$ and fixed $x$,
 \begin{equation}
 (\varepsilon_\alpha^{(\delta)}\ast\hat{u}_k)(x) \ra (\varepsilon_\alpha^{(\delta)}\ast\hat{u})(x).
 \label{d8}
 \end{equation}

Actually,
 \begin{eqnarray*}
 &&(\varepsilon_\alpha^{(\delta)}\ast\hat{u}_k)(x)- (\varepsilon_\alpha^{(\delta)}\ast\hat{u})(x)\\
&=&C\int_{|y-x|>\delta}\frac{\delta^\alpha[\hat{u}_k(y)-
\hat{u}(y)]}{(|x-y|^2-\delta^2)^{\frac{\alpha}{2}}|x-y|^n}dy\\
&=&C\{\int_{\substack{|y-x|>\delta\\ |y|<r-\eta}}
\frac{\delta^\alpha[\hat{u}_k(y)-
\hat{u}(y)]}{(|x-y|^2-\delta^2)^{\frac{\alpha}{2}}|x-y|^n}dy\\
&&+\underset{\substack{|y-x|>\delta \\r-\eta<|y|<r}}{\int}
\frac{\delta^\alpha[\hat{u}_k(y)-
\hat{u}(y)]}{(|x-y|^2-\delta^2)^{\frac{\alpha}{2}}|x-y|^n}dy\\
&&+
\int_{\substack{|y-x|>\delta \\|y|>r}}
\frac{\delta^\alpha[\hat{u}_k(y)-
\hat{u}(y)]}{(|x-y|^2-\delta^2)^{\frac{\alpha}{2}}|x-y|^n}dy\}\\
&=&C(I_1+I_2+I_3).
\end{eqnarray*}

For each fixed $x$ with $|x|<r$, choose $\delta$ and $\eta$ such that
$$B_\delta(x)\cap B^c_{r-2\eta}(0)=\emptyset.$$

It follows from (\ref{w1}) that as $k \ra \infty$
\begin{equation}
I_3=\int_{\substack{|y-x|>\delta \\
                      |y|>r}}\frac{\delta^\alpha[u_k(y)-
u(y)]}{(|x-y|^2-\delta^2)^{\frac{\alpha}{2}}|x-y|^n}dy \ra 0.
\end{equation}

\begin{eqnarray*}
I_2&=&\int_{\substack{ |y-x|>\delta \\
                      r-\eta<|y|<r}}
\frac{\delta^\alpha \int_{|z|>r}P_r(z,y)[u_k(z)-u(z)]dz}
{(|x-y|^2-\delta^2)^{\frac{\alpha}{2}}|x-y|^n}dy\\
&=&C\delta^\alpha\int_{|z|>r}\frac{u_k(z)-u(z)}{(|z|^2-r^2)^{\frac{\alpha}{2}}}
\int_{\substack{ |y-x|>\delta \\
                      r-\eta<|y|<r}}
\frac{(r^2-|y|^2)^{\frac{\alpha}{2}}dy}
{(|x-y|^2-\delta^2)^{\frac{\alpha}{2}}|x-y|^n|z-y|^n}dz\\
&=&C\delta^\alpha\int_{|z|>r}\frac{u_k(z)-u(z)}{(|z|^2-r^2)^{\frac{\alpha}{2}}}
\cdot I_{21}(x,z)dz.
\end{eqnarray*}
Noting that in the ring $r-\eta<|y|<r$, we have
$$|x-y|>\eta+\delta.$$
It then follows that
 \begin{eqnarray}
&&I_{21}(x,z)\nonumber\\
&\leq&\frac{1}{(2\eta\delta+\eta^2)^{\frac{\alpha}{2}}(\eta+\delta)^n}
\int_{r-\eta<|y|<r}\frac{(r^2-|y|^2)^{\frac{\alpha}{2}}dy}{|z-y|^n}
\nonumber\\
&=&C\int_{r-\eta}^r (r^2-\tau^2)^{\frac{\alpha}{2}}
\left\{\int_{S_\tau} \frac{1}{|z-y|^n} d\sigma_y\right\}d\tau\nonumber\\
&=&C\int_{r-\eta}^r (r^2-\tau^2)^{\frac{\alpha}{2}}
\left\{\int_0^\pi \frac{\omega_{n-2}(\tau\sin\theta)^{n-2}\tau d\theta}
{(\tau^2+|z|^2-2\tau|z|\cos\theta)^{\frac{n}{2}}}\right\}d\tau\nonumber\\
&=&C\int_{r-\eta}^r (r^2-\tau^2)^{\frac{\alpha}{2}}\frac{1}{\tau^n}
\int_0^\pi \frac{\tau^{n-1}\sin^{n-2}\theta d\theta}{(
(\frac{|z|}{\tau})^2-2\frac{|z|}{\tau}\cos\theta+1)^{\frac{n}{2}}}d\tau
\label{d5}\\
&=&C\int_{r-\eta}^r \frac{(r^2-\tau^2)^{\frac{\alpha}{2}}}{\tau}
\frac{d\tau}{(\frac{|z|}{\tau})^{n-2}((\frac{|z|}{\tau})^2-1)}
\int_{0}^\pi\sin^{n-2}\beta d\beta \label{d6}\\
&<&\frac{Cr^{n-1}}{|z|^{n-2}}\int_{r-\eta}^r
\frac{(r^2-\tau^2)^{\frac{\alpha}{2}}}{|z|^2-\tau^2} d\tau\nonumber\\
&=&\frac{Cr^{n-1}}{|z|^{n-2}}\cdot J.\nonumber
\end{eqnarray}
In the above, to derive (\ref{d6}) from (\ref{d5}), we have made the following substitution (See Appendix in \cite{L}):
$$\frac{\sin\theta}{\sqrt{(\frac{|z|}{\tau})^2-2\frac{|z|}{\tau}\cos\theta+1}}
=\frac{\sin\beta}{\frac{|z|}{\tau}},$$

To estimate the last integral $J$, we consider

\textit{(a)}\quad For $r<|z|<r+1$,
$$J\leq\int_{r-\eta}^r\frac{(r+\tau)^{\frac{\alpha}{2}-1}}{(r-
\tau)^{1-\frac{\alpha}{2}}}d\tau\leq C_{\alpha,r}.$$

\textit{(b)}\quad For $|z|\geq r+1$, obviously,
\begin{equation}
J\sim\frac{1}{|z|^2}, \mbox{ for $|z|$ large}.\nonumber
\end{equation}

In summary,
\begin{equation}
I_{21}(x,z)\sim\left\{\begin{array}{ll}
1,  &\mbox{ for $|z|$ near r},\\
|z|^n,&\mbox{ for $|z|$ large}.
\end{array}
\right.
\nonumber
\end{equation}

Therefore, by (\ref{w1}), as $k \ra \infty$,
\begin{equation}
I_2=\delta^\alpha\int_{|z|>r}\frac{u_k(z)-u(z)}{(|z|^2 -r^2)^{\alpha/2}} I_{21}(x,z)dz \ra 0.
\end{equation}

Now what remains is to estimate
$$I_1=\delta^\alpha\int_{|z|>r}\frac{u_k(z)-u(z)}{(|z|^2-r^2)^{\frac{\alpha}{2}}}
I_{11}(x,z)dz,$$
where
$$I_{11}(x,z)=\int_{\substack{|y-x|>\delta \\
                      |y|< r-\eta}}
 \frac{(r^2-|y|^2)^{\frac{\alpha}{2}}dy}
 {(|x-y|^2-\delta^2)^{\frac{\alpha}{2}}|x-y|^n|z-y|^n}.$$

 \begin{eqnarray}
I_{11}(x,z)
&\leq&\frac{r^\alpha}{\delta^n}\int_{\substack{ |y-x|>\delta \\
                      |y|< r-\eta}}
 \frac{dy}{(|x-y|^2-\delta^2)^{\frac{\alpha}{2}}|z-y|^n}\\
&\leq&\frac{r^\alpha}{\delta^n(|z|-r+\eta)^n}
\int_{\delta<|y-x|<2r}\frac{dy}
 {(|x-y|^2-\delta^2)^{\frac{\alpha}{2}}}\\
 &=&\frac{r^\alpha}{\delta^n(|z|-r+\eta)^n}
\int_\delta^{2r}
\frac{\omega_{n-1}\tau^{n-1}d\tau}{(\tau^2-\delta^2)^{\frac{\alpha}{2}}}\\
&\leq&\frac{C}{|z|^n}.
\end{eqnarray}

By (\ref{w1}), as $k \ra \infty$, we have $I_1 \ra 0.$ This verifies (\ref{d8}) and hence completes the proof.
\bigskip

\noindent \textbf{Proof of Lemma \ref{lem3.2}.}
By using the  property  $$\int_{|y-x|>r}\varepsilon^{(r)}_{\alpha}(x-y)=1,$$
we have
\begin{eqnarray}
&&\frac{1}{r^\alpha}\left[u(x)-\varepsilon^{(r)}_{\alpha}\ast u(x)\right]\nonumber\\
&=&\frac{1}{r^\alpha}u(x)- c\int_{|y-x|>r}\frac{u(y)}{(|x -y|^2-r^2)^{\frac{\alpha}{2}}|x-y|^n}dy\nonumber\\
&=&c\int_{|y-x|>r}
\frac{u(x)-u(y)}{(|x-y|^2-r^2)^{\frac{\alpha}{2}}|x-y|^n}dy.
\label{c1}
\end{eqnarray}

 Compare (\ref{c1}) with
$$(-\lap)^\frac{\alpha}{2}u(x)=\lim_{r\rightarrow0}\int_{|y-x|>r}
\frac{u(x)-u(y)}{|x-y|^{\alpha+n}}dy.$$
One may expect that
$$\lim_{r\rightarrow0}\int_{|y-x|>r}
\frac{u(x)-u(y)}{|x-y|^{\alpha+n}}dy=\lim_{r\rightarrow0}\int_{|y-x|>r}
\frac{u(x)-u(y)}{(|x-y|^2-r^2)^{\frac{\alpha}{2}}|x-y|^n}dy.$$

Indeed, consider
\begin{eqnarray}
&&\int_{|y-x|>r}
\frac{u(x)-u(y)}{|x-y|^{n}}\left(\frac{1}{(|x-y|^2-r^2)^{\frac{\alpha}{2}}}
-\frac{1}{|x-y|^{\alpha}}\right)dy\nonumber\\
&=&\int_{r<|y-x|<1}\frac{u(x)-
u(y)}{|x-y|^{n}}\left(\frac{1}{(|x-y|^2-r^2)^{\frac{\alpha}{2}}}
-\frac{1}{|x-y|^{\alpha}}\right)dy\nonumber\\
&&+\int_{|y-x|\geq1}\frac{u(x)-
u(y)}{|x-y|^{n}}\left(\frac{1}{(|x-y|^2-r^2)^{\frac{\alpha}{2}}}
-\frac{1}{|x-y|^{\alpha}}\right)dy\nonumber\\
&=&I_{1}+I_{2}.\label{c2}
\end{eqnarray}

It is easy to see that as $r\rightarrow0$, $I_2$ tends to zero. Actually, same conclusion is true for $I_1$.
\begin{eqnarray}
I_1&=&\int_{r<|y-x|<1}\frac{\nabla u(x)\cdot(y-x)+O(|y-x|^2)}{|x-
y|^{n}}\left(\frac{1}{(|x-y|^2-r^2)^{\frac{\alpha}{2}}}
-\frac{1}{|x-y|^{\alpha}}\right)dy\nonumber\\
\label{c2.5}\\
&\leq&C\int_{r<|y-x|<1}\frac{|x-
y|^2}{|x-y|^n}\left(\frac{1}{(|x-y|^2-r^2)^{\frac{\alpha}{2}}}
-\frac{1}{|x-y|^{\alpha}}\right)dy\label{c3}\\
&=&C\int_r^1\frac{\tau^2}{\tau^n}\left(\frac{1}{(\tau^
2-r^2)^{\frac{\alpha}{2}}}
-\frac{1}{\tau^{\alpha}}\right)\tau^{n-1}d\tau\label{c4}\\
&\leq&C\int_1^\infty\left(\frac{1}{r^\alpha(s^
2-1)^{\frac{\alpha}{2}}}
-\frac{1}{r^\alpha s^\alpha}\right)sr^2ds\label{c5}\\
&=&Cr^{2-\alpha}\int_1^\infty\left(\frac{s^\alpha-(s^
2-1)^{\frac{\alpha}{2}}}{(s^
2-1)^{\frac{\alpha}{2}}s^\alpha}\right)sds.\label{c6}
\end{eqnarray}
Equation (\ref{c2.5}) follows from the Taylor expansion. Due to symmetry, we have
$$\int_{r<|y-x|<1}\frac{\nabla u(x)\cdot(y-x)}{|x-
y|^{n}}\left(\frac{1}{(|x-y|^2-r^2)^{\frac{\alpha}{2}}}
-\frac{1}{|x-y|^{\alpha}}\right)dy=0$$
 and get (\ref{c3}). By letting $|y-x|=\tau$ and $\tau=rs$ respectively, one obtains (\ref{c4}) and (\ref{c5}).
It is easy to see that the integral in (\ref{c6}) converges near 1. To see that it also converges near infinity, we estimate
$$s^\alpha-(s^2-1)^{\frac{\alpha}{2}}.$$

Let$f(t)=t^{\alpha/2}$. By the \emph{mean value theorem},

\begin{eqnarray*}
f(s^2)-f(s^2-1)
&=&f^\prime(\xi)(s^2-(s^2-1))\\
&=&\frac{\alpha}{2}\xi^{\frac{\alpha}{2}-1}\sim s^{\alpha-2}, \mbox{ for $s$ sufficiently large.}
\end{eqnarray*}

This implies that
$$\frac{s^\alpha-(s^2-1)^{\frac{\alpha}{2}}}{(s^2-1)^{\frac{\alpha}{2}}s^\alpha}s
\sim\frac{s^{\alpha-2}s}{(s^2-1)^{\frac{\alpha}{2}}s^\alpha}
\sim\frac{1}{s^{1+\alpha}}.$$
Now it is obvious that (\ref{c6}) converges near infinity. Thus we have
$$\int_1^\infty\left(\frac{s^\alpha-(s^
2-1)^{\frac{\alpha}{2}}}{(s^
2-1)^{\frac{\alpha}{2}}s^\alpha}\right)sds<\infty.$$
Since $0<\alpha<2$, as $r\rightarrow0$, (\ref{c6}) goes to zero and
$I_1$ converges to zero. Together with (\ref{c1}) and (\ref{c2}), we get (\ref{c0}). This concludes the proof.
\bigskip

\noindent \textbf{Proof of Corollary \ref{corA}.}
Assume that $u \in L_{\alpha}$ is $\alpha$-harmonic (in the distribution sense).

From the Poisson's expression of $\hat{u}$, by the \emph{Lebesgue dominate convergence theorem}, one can differentiate under the integral signs to show that $\hat{u}$ is smooth in $B_r(0)$ for any $r>0$ (see also \cite{Fa}). Lemma \ref{lem3.1} indicates that $\hat{u}$ is $\alpha$-harmonic in the average sense, hence by Lemma \ref{lem3.2}, it is also $\alpha$-harmonic. That is we have
$$
\left\{\begin{array}{ll}
(-\Delta)^{\alpha/2} ( \hat{u} - u)(x) = 0 , & x \in B_r(0),\\
\hat{u} - u = 0 , & x \not{\in} B_r(0).
\end{array}
\right.
$$
Now by Silvestre's {\em maximum principle}, $u(x) = \hat{u}(x)$, hence $u$ is also $\alpha$-harmonic in the average sense.

The converse can be derive from Lemma \ref{lem3.2} immediately.

{\em Authors' Addresses and E-mails:}
\medskip

Wenxiong Chen

Department of Mathematical Sciences

Yeshiva University

New York, NY, 10033 USA

wchen@yu.edu
\medskip

Lorenzo D'Ambrosio

Dipartimento di Matematica

Universit\`a degli Studi di Bari

via Orabona, 4 - 70125 Bari - Italy

dambros@dm.uniba.it

\medskip

Yan Li

Department of Mathematical Sciences

Yeshiva University

New York, NY, 10033 USA

yali3@mail.yu.edu

\end{document}